\documentclass[a4paper,12pt]{amsart}

\usepackage{amsmath,amsfonts,amsthm,amsopn,cite,mathrsfs,amssymb}
\usepackage{mathtools}
\usepackage[T1]{fontenc}
\usepackage{lmodern}
\usepackage[a4paper,left=3cm,right=3cm,top=3.5cm,bottom=3.5cm]{geometry}
\usepackage{subfigure}

\newtheorem{theorem}{Theorem}[section]
\newtheorem{proposition}[theorem]{Proposition}
\newtheorem{lemma}[theorem]{Lemma}
\newtheorem{corollary}[theorem]{Corollary}
\newtheorem{definition}[theorem]{Definition}
\newtheorem{example}[theorem]{Example}
\theoremstyle{remark}
\newtheorem{remark}[theorem]{Remark}

\newcommand{\T}{\mathbb{T}}
\newcommand{\Z}{\mathbb{Z}}
\newcommand{\R}{\mathbb{R}}

\newcommand{\dd}{\,\mathrm{d}}
\newcommand{\e}{\mathrm{e}}

\begin{document}
\begin{abstract}
A closed rectifiable curve on the torus
$\mathbb T^d=\mathbb R^d/\mathbb Z^d$ is a $t$-design curve if
normalized arclength integration is exact for
trigonometric polynomials $e^{2\pi i \langle m,x\rangle}$ with frequencies $m\in\mathbb Z^d$
satisfying $\|m\|_\infty\le t$. 
For the arclength $L$ of such curves, we establish the sharp lower bound
\begin{equation}\label{eq:abs}
L\geq \Big(\frac{(t+1)^{2d}-1}{t(t+2)}\Big)^{1/2}.
\end{equation}
Equality is attained by the family of Korobov
curves described in \cite{Kunis,Lindblad}. 

For even $t$, we also derive an explicit lower bound on the length of
total-degree $t$-design curves, defined by the frequency condition
$\|m\|_1\leq t$. This bound is sharp for $t=2$ in every dimension
and for every even $t$ in dimension two.

We determine the minimum length of hyperbolic cross $t^2$-design
curves in dimension two. For $d\geq3$, we establish lower bounds
of optimal order for the length of hyperbolic cross $t^d$-design curves.

Finally, we obtain lower bounds of optimal order for the length of
$t$-design curves associated with sparse frequency sets.
\end{abstract}

\title[Optimal bounds for design curves on the torus]{Optimal bounds for design curves on the torus}
\author{Martin Ehler}
\address{Faculty of Mathematics \\
University of Vienna \\
Oskar-Morgenstern-Platz 1 \\
A-1090 Vienna, Austria}
\email{martin.ehler@univie.ac.at}
\maketitle

\vspace{-.5cm}

\section{Introduction}
Exact integration along closed curves was used in
\cite{EhlerGrafNeumayerSteidl} to approximate measures on the
torus and the sphere, with exactness imposed on trigonometric
polynomials and spherical polynomials, respectively. The notion
of a $t$-design curve was subsequently introduced on the sphere
in \cite{EhlerGrochenig}, requiring normalized arclength
integration to be exact for all spherical polynomials of
degree at most $t$. Such curves are analogues of spherical
$t$-design point sets \cite{DelsarteGoethalsSeidel}, which have
been studied extensively for several decades
\cite{BondarenkoRadchenkoViazovska}. More recently, asymptotically
optimal $t$-design curves were constructed by Lindblad
\cite{Lindblad} on $\mathbb S^3\subseteq\R^4$ and on tori,
and by Kunis \cite{Kunis} on tori. King and Mixon \cite{KingMixon} obtained explicit spectral lower bounds on the length and exact optimality results for $2$-design curves on odd-dimensional spheres $\mathbb S^{2d-1}\subseteq \R^{2d}$.  

We concentrate on the torus $\T^d=\R^d/\Z^d$.
\begin{definition}\label{def:1}
Let $t$ be a positive integer. A rectifiable curve \(\gamma:[0,L]\to \mathbb R^d\) with $\gamma(L)-\gamma(0)\in\mathbb Z^d$, parametrized by
arclength, is called a \emph{\(t\)-design curve} on $\T^d$ if
\begin{equation}\label{eq:design}
    \frac1L\int_0^L
    \e^{2\pi i \langle m,\gamma(s)\rangle } d s= \int_{\T^d} \e^{2\pi i \langle m,x\rangle }\dd x,
\end{equation}
for every \(m\in\Z^d\) satisfying $\|m\|_\infty\leq t$.
\end{definition}

Kunis \cite{Kunis} obtained explicit lower bounds on the length needed for exact trigonometric quadrature as in \eqref{eq:design}. The bounds are of optimal order in $t$ and allow additional weighting in the arclength integral, hence, hold for a larger class than \(t\)-design curves. For even $t$, a simplified version in 
\cite[Theorem~3.10]{Kunis} gives $L\geq \frac{0.94}{\sqrt{2d+1}}(t+1)^{d-1}$. Here, in the stricter unweighted setting of $t$-design curves, we derive the estimate \eqref{eq:abs} that implies $L\geq(t+1)^{d-1}$. This bound turns out to be sharp.

We adapt King and Mixon's spectral approach on the sphere \cite{KingMixon} to the torus:
\begin{itemize}
\item[(a)] The design property makes a suitable family of exponential functions along the curve orthonormal. 
\item[(b)] Fourier analysis and Bessel's inequality yield a lower bound for the averaged $L^2$-energy of derivatives of this family. 
\item[(c)] On the torus, the energy and the lower bound can be evaluated explicitly.
\end{itemize}

We derive results for four frequency shapes:
\subsubsection*{(1) $t$-design curve}
The resulting estimate \eqref{eq:abs} is sharp. The Korobov curves described in \cite{Kunis,Lindblad} attain this bound and therefore have minimum length among all
 $t$-design curves on $\mathbb T^d$. These main results are derived in Section \ref{sec:2}.  

\subsubsection*{(2) Total degree $t$-design curve}
We also consider the related notion of total-degree \(t\)-design curves in Section \ref{sec:total degree}, where the frequencies are $\{m\in\Z^d : \|m\|_1\leq t\}$. 
Curves of minimal length are constructed for total-degree \(1\) and \(2\) in every dimension and for every even total degree in dimension two. 

\subsubsection*{(3) Hyperbolic cross $t$-design curves} We discuss the hyperbolic cross~\cite{hyperbolic,hyperbolic2}, 
\begin{equation*}
H_d(t):=\{m\in\Z^d : \prod_{j=1}^d \max(1,|m_j|) \,\leq t\},
\end{equation*}
as an alternative to tensor-product frequency sets in Section \ref{sec:hyperbolic}. 
We derive sharp bounds for exact arclength integration of frequencies in the hyperbolic cross in dimension $d=2$. For $d\geq 3$, we obtain sharp lower bounds on the order of the length. This higher-dimensional hyperbolic-cross bound is attained by a construction that uses a uniform lattice-point counting argument and a hyperbolic-cross estimate.

\subsubsection*{(4) Sparse $t$-design curves}
Finally, we consider the sparse frequency sets
\[
    \Omega_{d,r}(t)
    =
    \left\{
        m\in\mathbb Z^d:
        \|m\|_\infty\le t,\quad
        \#\operatorname{supp}(m)\le r
    \right\},
    \qquad 1\le r\le d.
\]
The corresponding trigonometric polynomials are sums of terms
depending on at most $r$ coordinates. For fixed $d$ and $r$,
we establish lower bounds of order $t^{r-1}$ for the length of
design curves associated with these frequencies and verify the existence of 
simple closed geodesics with matching length order. 
Thus the optimal length exponent is determined by the sparsity
parameter $r$, rather than by the ambient dimension $d$.

\section{Sharp length bounds for $t$-design curves}\label{sec:2}
Let $d$ be a positive integer and denote the $d$-dimensional torus by $\T^d=\R^d/\Z^d$, endowed with the normalized Haar measure $\dd x$. We write $d s$ for the Lebesgue measure on the real line. A closed rectifiable curve on $\T^d$ of length $L>0$ is represented by a path $\gamma:[0,L]\to\mathbb{R}^d$, parametrized by
arclength, satisfying
\[
\gamma(L)-\gamma(0)\in\mathbb{Z}^d.
\]
The corresponding closed curve on $\T^d$ is obtained by
reduction modulo $\mathbb{Z}^d$.

The uniform condition $\|m\|_\infty\leq t$ follows the exact integration conditions used in \cite{Kunis}. Since the Fourier basis is orthonormal and the integral on the right-hand side of \eqref{eq:design} is the inner product with the constant $1$, we have $\int_{\T^d} \e^{2\pi i \langle m,x\rangle }\dd x=\delta_{m,0}$. 

\begin{theorem}\label{thm:lower bound}
Suppose $t$ is a positive integer. If $\gamma$ is a $t$-design curve on $\T^d$ of length $L$, then 
\begin{equation}\label{eq:LLL}
    L\ge \Big(\frac{(t+1)^{2d}-1}{t(t+2)}\Big)^{1/2} 
\end{equation}
There are $t$-design curves on $\T^d$ that attain this bound. 
\end{theorem}
The $t$-design curves attaining this bound appeared in \cite{Kunis,Lindblad}, and we refer to them as Korobov curves. 
Let $d$ and $t$ be positive integers, and define the Korobov vector
\begin{equation}\label{eq:Korobov}
v_{d,t}:=\big(1,t+1,\dots,(t+1)^{d-1}\big)^\top\in\mathbb Z^d.
\end{equation}
The associated curve
$\widetilde{\gamma}_{d,t}:[0,1]\to\mathbb R^d$, $\widetilde{\gamma}_{d,t}(s)= s\,v_{d,t}$ parametrizes the line segment from the origin to $v_{d,t}$. 
Its length is 
\begin{equation}\label{eq:length}
L=\|v_{d,t}\|_2 = \Big(\frac{(t+1)^{2d}-1}{t(t+2)}\Big)^{1/2}.
\end{equation}
The arclength parametrization is $\gamma_{d,t}:[0,L]\to\mathbb R^d$, $s\mapsto \frac{s}{\|v_{d,t}\|_2}v_{d,t}$, so that 
$\gamma_{d,t}(L)-\gamma_{d,t}(0)=v_{d,t}\in \mathbb{Z}^d$. It was mentioned in \cite{Lindblad} and verified in \cite{Kunis} that they form $t$-design curves. Thus, \eqref{eq:LLL} is sharp. 

The proof of the lower bound \eqref{eq:LLL} in Theorem \ref{thm:lower bound} needs some preparation. 
Let \(\gamma:[0,L]\to \mathbb R^d\) be a rectifiable curve, parametrized by
arclength. For every $a\in\mathbb R^d$, define a function $f_a:[0,L]\to\mathbb C$,
\begin{equation*}
f_a(s) := \e^{2\pi i \langle a,\gamma(s)\rangle}. 
\end{equation*}
We restrict $a$ to a suitable cube of side length $t$ centered at $0$. For every positive integer $t$, the set 
\begin{equation*}
[t]:=\{0,1,\dots,t\} - t/2
\end{equation*}
has cardinality $t+1$, and the Cartesian product $[t]^d$ has $(t+1)^d$ elements. If $t$ is even, then $[t]^d\subseteq \Z^d$, and $t$ odd yields $[t]^d\subseteq(\Z+\frac{1}{2})^d$. The difference satisfies $[t]^d-[t]^d\subseteq [-t,t]^d\cap\Z^d$, no matter if $t$ is even or odd. 

Invariance of $[t]^d$ under coordinate sign changes and permutations implies that the covariance matrix of the product lattice $[t]^d$ is a scalar multiple of the identity $I_d$. Further direct calculation for the discrete uniform distribution as in \cite[Section 2.3.3]{Freund} shows 
\begin{equation}\label{eq:11}
\frac{1}{(t+1)^d}\sum_{a\in[t]^d} a \,a^\top = \frac{t(t+2)}{12} I_d.
\end{equation}
This isotropy property is key to the following lemma.
\begin{lemma}\label{lemma:1}
Suppose $t$ is a positive integer and let \(\gamma:[0,L]\to\mathbb R^d\) be a rectifiable curve, parametrized by
arclength. Then $f_a'\in L^2([0,L],ds/L)$, for all $a\in [t]^d$, and 
\begin{equation}\label{eq:1}
\frac{1}{(t+1)^d} \sum_{a\in[t]^d} \frac{1}{L} \int_0^L |f_a'(s)|^2 d s =  4\pi^2  \frac{t(t+2)}{12}.
\end{equation}
\end{lemma}
The identity \eqref{eq:1} in Lemma \ref{lemma:1} holds independently of the design property but relies on the isotropy condition \eqref{eq:11}. 
\begin{proof}[Proof of Lemma \ref{lemma:1}]
Since $\gamma$ is rectifiable and parametrized by arclength, it is $1$-Lipschitz. Consequently, $f_a$ is Lipschitz as well. Therefore, both are differentiable almost everywhere with 
\begin{equation}\label{eq:arclength}
\|\gamma'\|_2=1,\qquad f_a' = 2\pi i \langle a,\gamma'\rangle f_a. 
\end{equation}
This also implies almost everywhere
\begin{equation}\label{eq:010}
|f_a'|^2 = 4\pi^2 |\langle a,\gamma'\rangle|^2 = 4\pi^2\,\gamma'^\top a \,a^\top \gamma'.
\end{equation}
We have $\|f_a\|_{L^\infty([0,L])} = 1$ and $\|f_a'\|_{L^\infty([0,L])}\leq 2\pi \|a\|_2$, so that, in particular, $f_a'\in L^2([0,L],\dd s/L)$. 

The relations \eqref{eq:010} and \eqref{eq:11} lead to
\begin{align*}
\frac{1}{(t+1)^d} \sum_{a\in[t]^d} \frac{1}{L} \int_0^L |f_a'(s)|^2 d s &= \frac{4\pi^2}{(t+1)^d L} \sum_{a\in[t]^d}  \int_0^L  \gamma'(s)^\top a\, a^\top \gamma'(s) \,d s\\
& = \frac{4\pi^2}{L}  \int_0^L  \gamma'(s)^\top \frac{t(t+2)}{12} I_d\, \gamma'(s) \,d s.
\end{align*}
Since $\gamma'(s)^\top \gamma'(s) = \|\gamma'(s)\|_2^2=1$ almost everywhere by \eqref{eq:arclength}, we arrive at \eqref{eq:1}. 
\end{proof}

We establish orthonormality of the functions $f_a$ for $a\in[t]^d$ as a consequence of the $t$-design property.

\begin{lemma}\label{lemma:2}
Let $t$ be a positive integer. If $\gamma$ is a $t$-design curve on $\T^d$ of length $L$, then the functions $\{f_a:a\in [t]^d\}$ are orthonormal in
$L^2([0,L],ds/L)$.
\end{lemma}
\begin{proof}[Proof of Lemma \ref{lemma:2}]
For $a,b\in[t]^d$, we have $a-b\in\Z^d$ and $\|a-b\|_\infty\leq t$, so that the $t$-design property \eqref{eq:design} yields
\begin{equation*}
\frac{1}{L} \int_0^L f_a(s)\overline{f_b(s)} d s  = \frac{1}{L} \int_0^L \e^{2\pi i \langle a-b,\gamma(s)\rangle} d s =\int_{\T^d} \e^{2\pi i \langle a-b,x\rangle} \dd x= \delta_{a,b}.
\end{equation*}
Thus, the functions $\{f_a : a\in[t]^d\}$ are orthonormal in $L^2([0,L],d s/L)$. 
\end{proof}
The orthonormality is key to deriving a lower bound on the left-hand side of \eqref{eq:1} for $t$-design curves. 
\begin{proposition}\label{prop:sobolev}
Let $t$ be a positive integer. If $\gamma$ is a $t$-design curve on $\T^d$ of length $L$, then
\begin{equation}\label{eq:2}
\frac{1}{(t+1)^d} \sum_{a\in[t]^d} \frac{1}{L} \int_0^L |f_a'(s)|^2 d s
     \ge \Big(\frac{2\pi}{L}\Big)^2 \frac{(t+1)^{2d}-1}{12}.
     \end{equation}
\end{proposition}
Theorem \ref{thm:lower bound} follows from combining Lemma \ref{lemma:1} with Proposition \ref{prop:sobolev}.
\begin{proof}[Proof of Theorem \ref{thm:lower bound}]
We estimate the left-hand side of \eqref{eq:1} by the lower bound in \eqref{eq:2} and obtain
\begin{equation*}
4\pi^2 \frac{t(t+2)}{12} \geq \Big(\frac{2\pi}{L}\Big)^2 \frac{(t+1)^{2d}-1}{12}. 
\end{equation*}
Canceling $4\pi^2/12$ and rearranging yields
\begin{equation*}
L^2\geq \frac{(t+1)^{2d}-1}{t(t+2)}. 
\end{equation*}
Taking square roots leads to \eqref{eq:LLL}. 
\end{proof}
To complete the proof of Theorem \ref{thm:lower bound}, we verify Proposition \ref{prop:sobolev}.

\begin{proof}[Proof of Proposition \ref{prop:sobolev}]
Let $N=\#[t]^d$ and $v=\gamma(L)-\gamma(0)\in \Z^d$. The parity of $t$ and $v_1+\dots+v_d$ lead to three cases. 

(i) We first treat the case when $t$ is even, so that $[t]^d\subseteq \mathbb Z^d$. For each $a\in[t]^d$, the function $f_a$ satisfies $f_a(L)=f_a(0)$, and we denote its Fourier coefficients by
\begin{equation*}
c_k(f_a) = \frac{1}{L}\int_0^L f_a(s) \e^{-2\pi i ks/L} d s,\qquad k\in\Z.
\end{equation*}
Since $f_a$ is Lipschitz continuous, it is absolutely continuous and the Fourier coefficients of $f'_a$ are determined by the standard relation 
\begin{equation}\label{eq:FT derive}
c_k(f'_a)=\frac{2\pi ik}{L} c_k(f_a),\qquad k\in\mathbb Z.
\end{equation}
Parseval's identity turns the $L^2$-norm of $f'_a$ into a weighted sequence norm of Fourier coefficients of $f_a$, 
\begin{equation*}
   \frac{1}{(t+1)^d} \sum_{a\in[t]^d} \frac1L\int_0^L\lvert f_a'(s)\rvert^2 d s
    = \frac{1}{(t+1)^d} \sum_{a\in[t]^d}\Big(\frac{2\pi}{L}\Big)^2\,\sum_{k\in\Z} k^2\lvert c_{k}(f_a)\rvert^2.
\end{equation*}

We define the total contribution of the \(k\)-th Fourier coefficient by
\[
    p_k:=\sum_{a\in[t]^d}\lvert c_{k}(f_a)\rvert^2.
\]

To bound $\sum_{k\in\Z} k^2 p_k$ from below, we use that the functions $\{f_a:a\in [t]^d\}$ are orthonormal by Lemma \ref{lemma:2}, so that Bessel's inequality yields
\begin{equation*}
p_k=\sum_{a\in[t]^d}\lvert c_{k}(f_a)\rvert^2\leq 1.
\end{equation*}
Parseval's identity for the Fourier basis implies 
\begin{equation*}
\sum_{k\in\Z}p_k = \sum_{k\in\Z} \sum_{a\in[t]^d}\lvert c_{k}(f_a)\rvert^2 = N.
\end{equation*}
Since $t$ is even, $N$ is odd, and $M:=(N-1)/2$ is an integer. We make use of 
\begin{equation*}
\sum_{|k|\leq M} (1-p_k)=\sum_{|k|>M}p_k,
\end{equation*}
to estimate $\sum_{k\in\Z} k^2 p_k$ from below, 
\begin{align*}
\sum_{k\in\Z} k^2 p_k - \sum_{|k|\leq M} k^2 & =
\sum_{|k| >M} k^2 p_k - \sum_{|k| \leq M} k^2 (1-p_k)\\
&\geq  \sum_{|k| >M} M^2 p_k - \sum_{|k| \leq M} M^2 (1-p_k) = 0.
\end{align*}
Consequently, $\sum_{k\in\Z} k^2 p_k  \geq \sum_{|k|\leq M} k^2$, which then implies the lower bound
\begin{equation}\label{eq:55}
\frac{1}{N}    \sum_{a\in[t]^d} \frac1L\int_0^L\lvert f_a'(s)\rvert^2 d s \geq \Big(\frac{2\pi}{L}\Big)^2 \frac{1}{N}\sum_{|k|\leq M} k^2.
\end{equation}
The sum evaluates to
\begin{equation*}
\frac{1}{N}\sum_{|k|\leq M} k^2 = \frac{2}{N}\sum_{k=1}^{(N-1)/2} k^2 = \frac{N^2-1}{12},
\end{equation*}
so that \eqref{eq:55} with $N=(t+1)^d$ leads to the claimed estimate \eqref{eq:2} in Proposition \ref{prop:sobolev}.

(ii) Let $t$ be odd and suppose that $v=\gamma(L)-\gamma(0)$ is such that $v_1+\dots+v_d$ is even. 
For $a\in[t]^d\subseteq (\mathbb Z +\frac{1}{2})^d$, this still implies $\langle a,v\rangle \in \Z$, as when $t$ is even. Hence,
\begin{align*}
f_a(L) = \e^{2\pi i \langle a,\gamma(L)\rangle} = \e^{2\pi i \langle a,v+\gamma(0)\rangle} = f_a(0).
\end{align*}
Thus, the functions $f_a$ still satisfy \eqref{eq:FT derive}, and we can follow the case (i) using the standard Fourier coefficients. 
Since $t$ is odd, $N$ is even, and $M:=N/2$ is an integer. To keep sums of length $N$, we switch from $\sum_{|k|\leq M}$ to $\sum_{k=-M}^{M-1}$. Computations analogous to the case (i) lead to 
\begin{equation*}
\sum_{k\in\Z} k^2 p_k  \geq \sum_{k=-M}^{M-1} k^2.
\end{equation*}
This yields the lower bound
\begin{equation*}
\frac{1}{N}    \sum_{a\in[t]^d} \frac1L\int_0^L\lvert f_a'(s)\rvert^2 d s \geq \Big(\frac{2\pi}{L}\Big)^2 \frac{1}{N}\sum_{k=-M}^{M-1}k^2.
\end{equation*}
The sum can be computed exactly by
\begin{equation*}
\frac{1}{N}\sum_{k=-M}^{M-1} k^2 =\frac{N^2+2}{12} \geq \frac{N^2-1}{12}.
\end{equation*}
This gives a slightly stronger bound in this case.

(iii) For $t$ odd and $v_1+\dots+v_d$ odd, we have $\langle a,v\rangle \in \Z+\frac{1}{2}$ for all $a\in[t]^d$, hence, 
the functions $f_a$ satisfy 
\begin{equation*}
f_a(L)= - f_a(0).
\end{equation*}
Thus, the derivative relation \eqref{eq:FT derive} for Fourier coefficients no longer holds with the standard Fourier basis, and we must change to an anti-periodic basis. We use
\begin{equation*}
\phi_k(s) = \e^{2\pi i (k+\frac{1}{2})s/L},\qquad k\in\Z.
\end{equation*}
This is an orthonormal basis for $L^2([0,L],\dd s/L)$, and the half-integer Fourier coefficients of $f_a$ are 
\begin{equation*}
c_{k+\frac{1}{2}}(f_a):=\frac{1}{L}\int_0^L f_a(s) \e^{-2\pi i (k+\frac{1}{2})s/L} d s,\qquad k\in\Z.
\end{equation*}
Since $f_a(L)=-f_a(0)$ and $e^{-2\pi i(k+\frac{1}{2})}=-1$, the boundary
term in integration by parts vanishes, and hence
\[
c_{k+\frac{1}{2}}(f_a')=\frac{2\pi i(k+\frac{1}{2})}{L}\,c_{k+\frac{1}{2}}(f_a).
\]
Thus, we obtain 
\begin{equation*}
   \frac{1}{(t+1)^d} \sum_{a\in[t]^d} \frac1L\int_0^L\lvert f_a'(s)\rvert^2 d s
    = \frac{1}{(t+1)^d} \sum_{a\in[t]^d}\Big(\frac{2\pi}{L}\Big)^2\,\sum_{k\in\Z} (k+\tfrac{1}{2})^2\lvert c_{k+\frac{1}{2}}(f_a)\rvert^2.
\end{equation*}
The total contribution of the $k+\frac{1}{2}$-th Fourier coefficient is 
\[
    p_{k+\frac{1}{2}}:=\sum_{a\in[t]^d}\lvert c_{k+\frac{1}{2}}(f_a)\rvert^2,\qquad k\in\Z.
\]
Bessel's inequality and Parseval's identity again give
\begin{equation*}
p_{k+\frac{1}{2}}=\sum_{a\in[t]^d}\lvert c_{k+\frac{1}{2}}(f_a)\rvert^2\leq 1,\qquad \sum_{k\in\Z}p_{k+\frac{1}{2}} 
= N. 
\end{equation*}
The same minimization argument as in case (i), now for half-integer frequencies and $M=N/2$, yields
\begin{equation*}
\frac{1}{N}    \sum_{a\in[t]^d} \frac1L\int_0^L\lvert f_a'(s)\rvert^2 d s \geq \Big(\frac{2\pi}{L}\Big)^2 \frac{1}{N}\sum_{k=-M}^{M-1}(k+\frac{1}{2})^2= \Big(\frac{2\pi}{L}\Big)^2 \frac{N^2-1}{12}.
\end{equation*}
This concludes the proof of Proposition \ref{prop:sobolev}.
\end{proof} 

\begin{remark}
Case (ii) in the proof of Proposition \ref{prop:sobolev} yields a stronger bound under an additional parity assumption. 
If $\gamma$ is a $t$-design curve on $\T^d$ of length $L$, $t$ is odd, and $v_1+\dots+v_d$ is even for $v=\gamma(L)-\gamma(0)$, then 
$ 
    L\ge  \big(\frac{(t+1)^{2d}+2}{t(t+2)}\big)^{1/2}$. 
\end{remark}

\section{Total-degree $t$-design curves}\label{sec:total degree}
Definition \ref{def:1} addressed the uniform frequency range $\|m\|_\infty\leq t$. Design curves are considered in \cite{Lindblad} for trigonometric polynomials of total degree at most \(t\), corresponding to the frequency condition $\|m\|_1\leq t$.
\begin{definition}
Let $t$ be a positive integer. A rectifiable curve \(\gamma:[0,L]\to \mathbb R^d\) with $\gamma(L)-\gamma(0)\in\mathbb Z^d$, parametrized by
arclength, is called a \emph{total-degree \(t\)-design curve} on $\T^d$ if
\begin{equation}\label{eq:design new}
    \frac1L\int_0^L
    \e^{2\pi i \langle m,\gamma(s)\rangle } d s= \int_{\T^d} \e^{2\pi i \langle m,x\rangle }\dd x,
\end{equation}
for every \(m\in\Z^d\) satisfying $\|m\|_1\leq t$.
\end{definition}
Every $t$-design curve in the sense of Definition \ref{def:1} is also a total-degree $t$-design curve, because $\|m\|_1\leq t$ implies $\|m\|_\infty\leq t$. 

The curve $\gamma_{d,t}$ induced by the Korobov vector \eqref{eq:Korobov} is a total-degree $(t+1)$-design curve. The latter is seen as follows. For $0<\|m\|_1\leq t+1$, either $\|m\|_\infty\leq t$, in which case exactness follows from the Korobov curve being a $t$-design, or there is $m_i=\pm (t+1)$, in which case $\langle m,v_{d,t}\rangle \neq 0$.

\subsection{Total-degree one}
We proceed by first determining the minimum length for total-degree \(1\)-design curves. For even degrees, we then derive a general lower bound and show that it is attained for \(t=2\) in every dimension and for every even \(t\) in dimension two.

\begin{lemma}\label{lemma:313}
The curve $\gamma:[0,\sqrt{d}]\to \R^d$, $s\mapsto\frac{s}{\sqrt{d}} ( 1,\dots,1)^\top$ is a total-degree $1$-design curve. Its length $\sqrt{d}$ is minimal among all total-degree $1$-design curves on $\T^d$. 
\end{lemma}
\begin{proof}
The design property follows from the observation that every nonzero $m\in\Z^d$ with $\|m\|_1\leq 1$ satisfies $\langle m,(1,\ldots,1)^\top\rangle\neq 0$. 

To prove the lower bound on the length, we again follow the approach
of \cite{KingMixon}, which takes a particularly simple form in this case. Let $\gamma:[0,L]\to \R^d$ be a total-degree $1$-design curve.  
For $j=1,\dots,d$, define $f_j:=f_{e_j}$, so that 
\begin{equation*}
f_j(s) = \e^{2\pi i \langle e_j,\gamma(s)\rangle} = \e^{2\pi i \gamma_j(s)} ,\qquad 0\leq s\leq L.
\end{equation*} 
According to $|f'_j|^2 = 4\pi^2 |\gamma_j'|^2$ and $1=\|\gamma'\|_2^2=\sum_{j=1}^d |\gamma'_j|^2 $ almost everywhere we obtain
\begin{equation*}
\frac{1}{d}\sum_{j=1}^d\frac{1}{L}\int_0^L | f_j'(s)|^2 d s  = \frac{4\pi^2}{d}.
\end{equation*}

Since $|f_j(s)|=1$, the functions $f_j$ have unit norm, $\frac{1}{L}\int_0^L |f_j(s)|^2\dd s = 1$. The design property does not imply orthonormality of $f_1,\dots f_d$, but they are orthogonal to the constant function, $\frac{1}{L}\int_0^L f_j(s)\dd s = 0$. 
A Fourier coefficient argument using \eqref{eq:FT derive} leads to 
\begin{equation*}
\frac{1}{L}\int_0^L | f_j'(s)|^2 \dd s \geq \Big(\frac{2\pi}{L}\Big)^2.
\end{equation*}
Thus, we obtain $\frac{4\pi^2}{d} \geq (\frac{2\pi}{L})^2$, which yields $L\geq \sqrt{d}$.
\end{proof}
\subsection{Even total-degree}
Assume $t=2n$ is even. Set
\begin{align*}
A_{d,n}:=	\{a\in\Z^d : \|a\|_1\leq n\},\qquad N_{d,n}:=\# A_{d,n}.	
\end{align*}
This yields $A_{d,n}-A_{d,n}\subseteq \{m\in\Z^d : \|m\|_1\leq t\}$, and 
the quantity 
\begin{equation*}
\sigma_{d,n}:=\frac{1}{N_{d,n}} \sum_{a\in A_{d,n}} a_j^2,\qquad j=1,\dots,d,
\end{equation*}
is independent of $j$. Invariance of \(A_{d,n}\) under coordinate sign changes and permutations implies
\begin{equation*}
\frac{1}{N_{d,n}} \sum_{a\in A_{d,n}} a\,a^\top = \sigma_{d,n} \,I_d,
\end{equation*}
that leads to a lower bound on the length of total-degree $t$-design curves.

\begin{theorem}\label{thm:second}
Suppose $t=2n$ is a positive even integer. If $\gamma$ is a total-degree $t$-design curve on $\T^d$ of length $L$, then 
\begin{equation}\label{eq:LLLnew}
    L\ge  \Big(\frac{N_{d,n}^{2}-1}{12 \sigma_{d,n}}\Big)^{1/2}.
\end{equation}
\end{theorem}
\begin{proof}
We repeat the argument of Section \ref{sec:2} with the functions $f_a(s) = \e^{2\pi i \langle a,\gamma(s)\rangle}$ but the centered cube $[t]^d$ replaced by $A_{d,n}$. The isotropy condition allows the computation of 
\begin{equation*}
\frac{1}{N_{d,n}} \sum_{a\in A_{d,n}} \frac{1}{L} \int_0^L |f_a'(s)|^2  d s = 4\pi^2 \sigma_{d,n},
\end{equation*}
and the total-degree $t$-design property \eqref{eq:design new} ensures that the functions $\{f_a : a\in A_{d,n}\}$ are orthonormal in $L^2([0,L],d s/L)$.  

The set \(A_{d,n}\) consists of the origin and pairs \(\{a,-a\}\), so that \(N_{d,n}\) is odd, and we may follow case (i) in the proof of Proposition \ref{prop:sobolev}. It leads to 
\begin{equation*}
\frac{1}{N_{d,n}} \sum_{a\in A_{d,n}} \frac{1}{L} \int_0^L |f_a'(s)|^2 d s \geq \Big(\frac{2\pi}{L}\Big)^2 \frac{N_{d,n}^2-1}{12}.
\end{equation*}
Thus, we arrive at 
\begin{equation*}
4\pi^2 \sigma_{d,n} \geq \Big(\frac{2\pi}{L}\Big)^2\frac{N_{d,n}^2-1}{12},
\end{equation*}
which leads to the estimate \eqref{eq:LLLnew}. 
\end{proof}

The quantities $N_{d,n}=\# A_{d,n}$ are known as Delannoy numbers \cite{Delannoy},  
\begin{align*}
N_{d,n} & = \sum_{j=0}^{\min(d,n)} 2^j \binom{d}{j} \binom{n}{j}.
\end{align*}
This is a polynomial of degree $d$ in $n$ with leading coefficient $2^d/d!$. Hence, for large $n$, the asymptotics $N_{d,n}\sim 2^dn^d/d!$ hold.
\begin{lemma}
For every \(d,n\ge1\), we have $\sigma_{d,n}= \frac{1}{2}\big(\frac{N_{d+2,n}}{N_{d,n}}-1 \big)$, and $n\to\infty$ yields the leading term $\sigma_{d,n}\sim \frac{2}{(d+1)(d+2)}n^2$, so that
\begin{equation}\label{eq:asym total}
\Big(\frac{N_{d,n}^{2}-1}{12 \sigma_{d,n}}\Big)^{1/2}\sim \frac{2^d}{d!}\,\Big(\frac{(d+1)(d+2)}{24}\Big)^{1/2}n^{d-1}.
\end{equation}
\end{lemma}
\begin{proof}
Fixing the first coordinate at $a_1=r$ in $\sigma_{d,n}=\frac{1}{N_{d,n}} \sum_{a\in A_{d,n}} a_1^2$ leaves exactly $N_{d-1,n-|r|}$ choices for the remaining coordinates, so that 
\begin{equation*}
    \sigma_{d,n}
    = \frac{1}{N_{d,n}}
      \sum_{r=-n}^{n} r^2 N_{d-1,n-|r|}
    = \frac{2}{N_{d,n}}
      \sum_{r=1}^{n} r^2 N_{d-1,n-r},
\end{equation*}
where we use the convention $N_{0,p}=N_{q,0}=1$ for $p,q\geq 0$. 

For $r\geq 1$, the number of triplets $(a_1,a_2,a_3)\in\Z^3$ satisfying $|a_1|+|a_2|+|a_3|=r$ can be counted by adding the numbers for exactly one, two, or three nonzero coordinates and is $6+12(r-1)+8\binom{r-1}{2}=4r^2+2$. According to the sum \(r\) of the absolute values of their last three coordinates, we obtain
\begin{equation*}
N_{d+2,n} = N_{d-1,n}+\sum_{r=1}^n(4r^2+2)N_{d-1,n-r}.
\end{equation*}
Similarly, we determine
\begin{equation*}
N_{d,n} = N_{d-1,n}+2\sum_{r=1}^n N_{d-1,n-r}.
\end{equation*}
Subtracting one from the other yields 
\begin{equation*}
N_{d+2,n} - N_{d,n} = 4\sum_{r=1}^n r^2 N_{d-1,n-r}=2N_{d,n}\sigma_{d,n},
\end{equation*}
so that $\sigma_{d,n}=\frac{1}{2}\big(\frac{N_{d+2,n}}{N_{d,n}}-1 \big)$. 

The observation $N_{d,n}\sim 2^dn^d/d!$ implies $\sigma_{d,n}\sim \frac{2}{(d+1)(d+2)}n^2$ and combining both gives \eqref{eq:asym total}. 
\end{proof}

\subsection{Minimal length for $d=2$ or $t=2$}
We specify Theorem \ref{thm:second} for two distinct cases. 
\begin{corollary}\label{cor:d=2}
Suppose $n\geq 1$ is an integer. 
\begin{itemize}
\item[($d=2$)] If $\gamma$ is a total-degree $2n$-design curve on $\T^2$ of length $L$, then 
\begin{equation}\label{eq:LLLnew2}
    L\ge  \sqrt{2n^2+2n+1}.
\end{equation}

\item[($t=2$)] If $\gamma$ is a total-degree $2$-design curve on $\T^d$ of length $L$, then 
\begin{equation}\label{eq:LLLnew3}
    L\ge  \sqrt{\tfrac{1}{6}\,d(d+1)(2d+1)}.
\end{equation}
\end{itemize}
\end{corollary}
We will see that there are curves that attain these bounds.
\begin{proof}
(a) For $d=2$, we have $N_{2,n}=2n^2+2n+1$ and $N_{4,n}-N_{2,n}=\frac{2}{3}n(n+1)(n^2+n+1)$. We can calculate $\sigma_{2,n}$ by 
\[
    \sigma_{2,n}
    =\frac{N_{4,n}-N_{2,n}}{2N_{2,n}}
    =  \frac{n(n+1)(n^2+n+1)}
         {3(2n^2+2n+1)}. 
\]
Since $N_{2,n}^2-1=4n(n+1)(n^2+n+1)$, substitution into \eqref{eq:LLLnew} gives $L^2\geq \frac{N_{2,n}^{2}-1}{12 \sigma_{2,n}}= 2n^2+2n+1$.

(b) For $t=2$, the Delannoy numbers simplify to $N_{d,1}=2d+1$, and we obtain $\sigma_{d,1}=2/(2d+1)$, so that $\frac{N_{d,1}^{2}-1}{12 \sigma_{d,1}} = \frac{1}{6}d(d+1)(2d+1)$. 
\end{proof}

\begin{example}\label{ex:ex 1}
Let $t=2$. The curve associated to the vector 
\begin{equation*}
w_{d} = (1,2,\dots,d)^\top 
\end{equation*}
is a total-degree \(2\)-design curve on $\T^d$. Indeed, if \(m\) has one nonzero entry, then \(\langle m,w_d\rangle\ne0\). If it has two nonzero entries, both have modulus one, and the corresponding entries of \(w_d\) are positive and distinct, giving \(\langle m,w_d\rangle\ne0\). Its length is  
$L_d=\sqrt{\frac{1}{6}\,d(d+1)(2d+1)}$, 
which is minimal among all total-degree $2$-design curves on $\T^d$ by Corollary \ref{cor:d=2}.
\end{example}

The following construction provides total-degree \(2n\)-design curves in every dimension and attains the bound \eqref{eq:LLLnew2} of Corollary \ref{cor:d=2} when $d=2$.

Let $t$ be a positive even integer, so that $t=2n$. The curve associated to the vector 
\begin{equation*}
w_{d,n} = \big( n^{d-1},n^{d-2} (n+1),\dots, (n+1)^{d-1}\big)^\top 
\end{equation*}
is $\beta_{d,n}(s)=\frac{s}{\|w_{d,n}\|_2}w_{d,n}$, for $0\leq s\leq \|w_{d,n}\|_2=L_{d,n}$. 
\begin{lemma}\label{lemma:betadn}
The curve $\beta_{d,n}$ is a total-degree \(2n\)-design curve on $\T^d$ of length 
\begin{equation*}
L_{d,n}=\Big(\frac{(n+1)^{2d}-n^{2d}}{2n+1}\Big)^{1/2}.
\end{equation*}
\end{lemma}
For large $n$, the asymptotic is $L_{d,n}\sim \sqrt{d}\,n^{d-1}$.
\begin{proof}
Suppose $\langle m,w_{d,n}\rangle=0$ for some nonzero $m\in\mathbb{Z}^d$.
Since every entry of $w_{d,n}$ is positive, $m$ has at least two nonzero
entries. Let $r<s$ be their first and last indices. Dividing the
relation by $n^{d-s}(n+1)^{r-1}$ gives
\[
\sum_{j=r}^{s} m_j n^{s-j}(n+1)^{j-r}=0.
\]
Every term with $j>r$ contains a factor $n+1$. Since the sum is zero,
$n+1$ must also divide the first term $m_r n^{s-r}$. As $n$ and $n+1$
are coprime, this gives $(n+1)\mid m_r$. Similarly, every term with
$j<s$ contains a factor $n$, so $n$ divides the last term
$m_s(n+1)^{s-r}$, and hence $n\mid m_s$. Both entries
are nonzero, so 
\begin{equation*}
\|m\|_1\ge |m_r|+|m_s|\ge 2n+1.
\end{equation*}
Thus, whenever $0<\|m\|_1\le 2n$, the integer $\langle m,w_{d,n}\rangle$
is nonzero. It follows that $\beta_{d,n}$ is a total-degree \(2n\)-design curve on $\T^d$. 

To compute the length $L_{d,n}=\|w_{d,n}\|_2$, we recall the standard identity
\begin{equation*}
(a-b)\sum_{j=0}^{d-1} a^{d-1-j} b^j = a^d-b^d,\qquad a\neq b.
\end{equation*}
For $a=(n+1)^2$ and $b=n^2$, this yields 
\begin{equation*}
\|w_{d,n}\|_2^2 = \sum_{j=0}^{d-1} (n+1)^{2(d-1-j)}n^{2j} = \frac{(n+1)^{2d}-n^{2d}}{2n+1}.\qedhere
\end{equation*}
\end{proof}

\begin{example}\label{ex:opt irgendwie}
For $d=2$ and $n$ a positive integer, $\beta_{2,n}$ is a total-degree \(2n\)-design curve on $\T^2$ with $w_{2,n}=(n,n+1)^\top$ and length $L_{2,n}=\sqrt{2n^2+2n+1}$. The length is minimal among all total-degree $2n$-design curves on $\T^2$ by Corollary \ref{cor:d=2}.
\end{example}

Examples \ref{ex:ex 1} and \ref{ex:opt irgendwie} show that the bounds in Corollary \ref{cor:d=2} are sharp. 

\section{Hyperbolic cross design curves}\label{sec:hyperbolic}

Hyperbolic crosses provide a standard alternative to tensor-product
frequency sets \cite{hyperbolic,hyperbolic2}. We denote the $d$-dimensional hyperbolic cross by
\begin{equation*}
H_d(t):=\{m\in\Z^d : \prod_{j=1}^d \max(1,|m_j|) \,\leq t\},
\end{equation*}
and define design curves for frequencies in $H_d(t)$.
\begin{definition}
Let $t$ be a positive integer. A rectifiable curve \(\gamma:[0,L]\to \mathbb R^d\) with $\gamma(L)-\gamma(0)\in\mathbb Z^d$, parametrized by
arclength, is called a \emph{hyperbolic cross \(t\)-design curve} on $\T^d$ if
\begin{equation*}
    \frac1L\int_0^L
    \e^{2\pi i \langle m,\gamma(s)\rangle } d s= \int_{\T^d} \e^{2\pi i \langle m,x\rangle }\dd x,
\end{equation*}
for every \(m\in H_d(t)$.
\end{definition}
\subsection{Sharp bounds for $d=2$}
For the hyperbolic cross in dimension two, we obtain a sharp lower bound on the length of such design curves. 
\begin{theorem}
\label{thm:hyperbolic-cross-d2}
Let $n\geq 1$ be an integer and $d=2$.
\begin{itemize}
\item[(a)] If $\gamma$ is a hyperbolic cross $n^2$-design curve on $\T^2$, then its length $L$ satisfies 
\[
    L\geq \sqrt{2n^2+2n+1}.
\]
\item[(b)] The curve $\beta_{2,n}$ in Example \ref{ex:opt irgendwie} is a hyperbolic cross $n^2$-design and attains this bound.
\end{itemize}
\end{theorem}
It is a feature of the hyperbolic cross that the design strength can be $n^2$ although the length is only of order $n$. 
\begin{proof}
(a) Suppose first that $n\geq2$. Consider $m\in\mathbb{Z}^2$ such that $\|m\|_1\leq2n$. If both coordinates are nonzero, then \(|m_1m_2|\leq(\|m\|_1/2)^2\leq n^2\). Otherwise we have 
\begin{equation*}
\max\{1,|m_1|\}\max\{1,|m_2|\}\leq 2n\leq n^2.
\end{equation*}
Consequently,
\[
    \{m\in\mathbb{Z}^2:\|m\|_1\leq2n\}
    \subseteq H_2(n^2),
\]
and every hyperbolic cross
$n^2$-design curve on $\T^2$ is also a total-degree $2n$-design
curve. Corollary \ref{cor:d=2} therefore yields $L\geq\sqrt{2n^2+2n+1}$.

For $n=1$, the hyperbolic cross is $H_2(1)=\{-1,0,1\}^2$, so $\gamma$ is a $1$-design curve in the sense of
Definition \ref{def:1}. Theorem \ref{thm:lower bound} gives $L\geq
  (\frac{2^4-1}{3})^{1/2}
    =\sqrt{5}$, which is the claimed bound for $n=1$.

(b) The curve in Example~\ref{ex:opt irgendwie} is $\beta_{2,n}(s)=\frac{s}{L_{2,n}}w_{2,n}$ with $w_{2,n}=(n,n+1)^\top$ and its length $L_{2,n}=\sqrt{2n^2+2n+1}$ matches the bound in Part (a). To verify the design property, it is sufficient to check that 
\[
    \langle m,w_{2,n}\rangle\neq0
    \qquad\text{for every }m\in H_2(n^2)\setminus\{0\}.
\]
Indeed, if $m\neq 0$ and $nm_1+(n+1)m_2=0$, then $(n+1)m_2=-nm_1$. Since $n$ and
$n+1$ are coprime, $n+1$ must divide $m_1$, so that there is $k\in \Z\setminus\{0\}$ such that $k(n+1)=m_1$. Thus, we obtain 
\begin{equation*}
0 = nk(n+1)+(n+1)m_2 = (n+1)(nk+m_2),
\end{equation*}
so that $m_2=-nk$. This also implies $m = k(n+1,-n)^\top$. Consequently,
\[
    \max\{1,|m_1|\}\max\{1,|m_2|\}
    =k^2n(n+1)
    \geq n(n+1)>n^2,
\]
which leads to $m\not \in H_2(n^2)$.

Thus, for every nonzero $m\in H_2(n^2)$, the integer
$\langle m,w_{2,n}\rangle$ is nonzero, and this implies exact integration of $e^{2\pi i \langle m,x\rangle}$. For $m=0$, the normalized integral equals $1$. Hence
$\beta_{2,n}$ is a hyperbolic cross $n^2$-design curve. 
\end{proof}
\subsection{Sharp asymptotic bounds for $d\geq 3$}
We determine a sharp asymptotic bound on the length when the dimension is bigger than two.
\begin{theorem}\label{thm:hyperbolic asymptotic}
Let $d\geq 3$ be an integer. 
\begin{itemize}
\item[(a)] If $\gamma$ is a hyperbolic cross $n^d$-design curve on $\T^d$, for some positive integer $n$, then its length $L$ satisfies 
\[
    L\geq n^{d-1}.
\]
\item[(b)] There is a constant $C_d>0$ such that, for all positive integers $n$, there exists a simple, hyperbolic cross $n^d$-design curve on $\T^d$ of length $L\leq C_d n^{d-1}$.
\end{itemize}
\end{theorem}

Before we can prove this theorem, we need to verify some estimates regarding the hyperbolic cross.
\begin{lemma}\label{lem:hc-reciprocal-sum}
Let $n\ge1$ be an integer.
\begin{itemize}
\item[(i)] For all integers $d\geq 2$,
\begin{equation}\label{eq:crazy inequality}
    \sum_{m\in H_d(n^d)\setminus\{0\}}
        \frac1{\|m\|_\infty}
    \leq \frac{(3d)^d}{d-1}\, n^{d-1}.
\end{equation}
\item[(ii)] For all integers $d\geq 3$, there is a constant $c_d>0$ such that 
\begin{equation}\label{eq:cardinality hc}
  \#H_d(n^d)  \leq c_d  n^{2(d-1)}.
\end{equation}
\end{itemize}

\end{lemma}

\begin{proof}[Proof of Lemma \ref{lem:hc-reciprocal-sum}]
(i) Let $N=n^d$. For $m\in H_d(N)\setminus\{0\}$, the substitution $a_j=\max(1,|m_j|)$ has at most $3^d$
preimages and satisfies $\|a\|_\infty=\|m\|_\infty$, so that we derive
\[
  \sum_{m\in H_d(N)\setminus\{0\}}
        \frac1{\|m\|_\infty} \leq 3^d 
    \sum_{\substack{a_1,\dots,a_d\geq 1\\
                     a_1\cdots a_d\le N}}
        \frac1{\|a\|_\infty}.
\]
For every summation index $a$, we define the boxes $Q_a:=\prod_{j=1}^d(a_j-1,a_j]$, so that $x\in Q_a$ satisfies
\[
    x_1\cdots x_d\le N,
    \qquad
    \frac1{\|x\|_\infty}\ge\frac1{\|a\|_\infty}.
\]
Since these boxes are disjoint and have volume one,
\[
  \sum_{\substack{a_1,\dots,a_d\geq 1\\
                     a_1\cdots a_d\le N}}
        \frac1{\|a\|_\infty}\le
    \int_{\substack{x\in(0,\infty)^d\\
                     x_1\cdots x_d\le N}}
        \frac{dx}{\|x\|_\infty}.
\]
Except on a set of measure zero, every \(x\) has a unique largest coordinate. Thus we can split the integration domain into \(d\) regions, according to which coordinate is largest. By symmetry, we may take $x_d=r$ as the largest coordinate
and write $x_j=ry_j$, with $0<y_j<1$, for $j<d$. 
According to $dx = r^{d-1}dy dr$ and $x_d=r$, we obtain
\begin{align*}
   \int_{\substack{x\in(0,\infty)^d\\
                     x_1\cdots x_d\le N}}
        \frac{dx}{\|x\|_\infty}
        & = d\int_{\substack{x\in(0,\infty)^d\\
                     x_1\cdots x_d\le N\\ x_1,\dots,x_{d-1}<x_d}}
        \frac{dx}{x_d}
        \\
    &=
    d\int_{(0,1)^{d-1}}
      \int_0^{(N/(y_1\cdots y_{d-1}))^{1/d}}
          r^{d-2}\,dr\,dy\\
          &=
    \frac{d}{d-1}N^{(d-1)/d}
    \int_{(0,1)^{d-1}} (y_1\cdots y_{d-1})^{-(d-1)/d}\,dy\\
    &=
    \frac{d}{d-1}N^{(d-1)/d}
    \left(\int_0^1 s^{-(d-1)/d}\,ds\right)^{d-1}
    \intertext{Since $\int_0^1 s^{-(d-1)/d}\,ds=\big[d\,s^{1/d}\big]_0^1=d$, we derive}
     \int_{\substack{x\in(0,\infty)^d\\
                     x_1\cdots x_d\le N}}
        \frac{dx}{\|x\|_\infty}&=
    \frac{d^d}{d-1}N^{(d-1)/d}.
\end{align*}
Multiplication by $3^d$ yields \eqref{eq:crazy inequality}. 

(ii) All implicit constants\footnote{We write \(\lesssim_d\) to indicate that the left-hand side is bounded above by a constant multiple of the right-hand side, where the constant may depend on \(d\).} below depend only on $d$. 
According to \cite{hyperbolic}, we note that 
\begin{equation*}
\#H_d(N)\lesssim_d N(1+\log N)^{d-1}
\end{equation*}
provided that $N\ge1$ and $d\geq 3$. 
Since $(1+\log N)^{d-1}\lesssim_d N^{(d-2)/d}$, we obtain
\begin{equation*}
  \#H_d(N)  \lesssim_d N^{1+(d-2)/d} = N^{2(d-1)/d}.\qedhere
\end{equation*}
\end{proof}

\begin{proof}[Proof of Theorem \ref{thm:hyperbolic asymptotic}]
(a) We observe that $[-n,n]^d\cap\Z^d\subseteq H_d(n^d)$, so that Theorem \ref{thm:lower bound} yields 
\begin{align*}
L^2& \ge \frac{(n+1)^{2d}-1}{n(n+2)} =\sum_{j=0}^{d-1}(n+1)^{2j} \geq (n+1)^{2(d-1)}.
\end{align*}
This implies $L\geq (n+1)^{d-1}\geq n^{d-1}$.

(b) We aim to find $v\in\Z^d$ such that $\|v\|_2\leq C_d n^{d-1}$ and $\langle v,m\rangle\neq 0$, for all $m\in H_d(n^d)\setminus\{0\}$. Then the straight line $s\mapsto \frac{s}{\|v\|_2} v$ constitutes the hyperbolic cross $n^d$-design curve.

We say that $m\in\mathbb Z^d\setminus\{0\}$ is primitive if its coordinates have no common divisor greater than $1$, meaning $\mathop{gcd}(|m_1|,\dots,|m_d|)=1$. Let $\mathcal P_d(n^d)$ be the primitive vectors in
$H_d(n^d)\setminus\{0\}$. Every forbidden hyperplane is represented
by a vector in $\mathcal P_d(n^d)$, so we may define, for every integer $k\ge1$, 
\begin{equation*}
\mathcal{B}_{k,m}:=\{v\in[-k,k]^d\cap\mathbb Z^d:\langle v,m\rangle=0\},\qquad m\in \mathcal P_d(n^d).
\end{equation*}
The complete forbidden set is $\mathcal{B}_k:=\bigcup_{m\in \mathcal P_d(n^d)} \mathcal{B}_{k,m}$.

The hyperplane-counting bound \cite[Equation (14)]{Forst} is
\begin{equation}\label{eq:hyperplane-count}
    \#\mathcal{B}_{k,m}
    \le
    \Big(\frac{2k}{\|m\|_\infty}+1\Big)(2k+1)^{d-2}.
\end{equation}
The union bound $\# \mathcal{B}_k \leq \sum_{m\in \mathcal P_d(n^d)}\# \mathcal{B}_{k,m}$ leads to
\begin{align*}
\# \mathcal{B}_k &\leq \sum_{m\in \mathcal P_d(n^d)} \Big(\frac{2k}{\|m\|_\infty}+1\Big)(2k+1)^{d-2}\\
&= (2k+1)^{d-2} \Big(\# \mathcal P_d(n^d) + \sum_{m\in \mathcal P_d(n^d)} \frac{2k}{\|m\|_\infty} \Big).
\end{align*}
The cardinality of the integer cube $[-k,k]^d\cap\mathbb Z^d$ is $(2k+1)^d$, so that the proportion of forbidden vectors in the cube is
\begin{equation*}
\frac{\# \mathcal{B}_k }{(2k+1)^d} \leq  \frac{\# \mathcal P_d(n^d) + 2k\sum_{m\in \mathcal P_d(n^d)} \frac{1}{\|m\|_\infty} }{(2k+1)^2}.
\end{equation*}
Lemma \ref{lem:hc-reciprocal-sum} with $\mathcal P_d(n^d)\subseteq H_d(n^d)$ implies the estimates
\begin{align*}
\frac{\# \mathcal{B}_k }{(2k+1)^d} &\lesssim_d  \frac{ n^{2(d-1)} + 2k n^{d-1}}{(2k+1)^2}\\
& \leq \Big(\frac{n^{d-1}}{k}\Big)^2 +\frac{n^{d-1}}{k}.
\end{align*}

Thus, there is a constant $A_d>0$ such that this proportion is smaller than one for $k=\lceil A_d n^{d-1}\rceil$. Hence, there is 
a vector $w\in[-k,k]^d\cap \Z^d$ with $w\not\in \mathcal{B}_k$. This implies that $\langle w,m\rangle\neq 0$ for all $m\in \mathcal P_d(n^d)$ and therefore also for all $m\in H_d(n^d)\setminus\{0\}$.  Its norm satisfies
\[
    \|w\|_2\le\sqrt d\,\|w\|_\infty \leq  \sqrt d k
    \le\sqrt d\,(A_d+1)n^{d-1}.
\]
We may put $C_d=\sqrt d\,(A_d+1)$. 
To obtain a simple curve, we use 
\begin{equation*}
v:=w/\mathop{gcd}(|w_1|,\dots,|w_d|). 
\end{equation*}
The length bound $L:=\|v\|_2\leq C_d n^{d-1}$ still holds and so do the nonorthogonality conditions. Thus, the curve $\gamma:[0,L]\to\T^d$,  $s\mapsto \frac{s}{\|v\|_2}v$ is simple and a hyperbolic $n^d$-design curve on $\T^d$. 
\end{proof}

\section{Sparse frequency sets}\label{sec:sparse}
For integers $d\ge2$, $1\le r\le d$, and $t\ge1$, define the
sparse frequency set
\[
    \Omega_{d,r}(t)
    :=
    \left\{
        m\in\mathbb Z^d:
        \|m\|_\infty\le t,\quad
        \#\operatorname{supp}(m)\le r
    \right\},
\]
where $\operatorname{supp}(m):=\{j:m_j\ne0\}$.

\begin{definition}\label{def:sparse-design}
For integers $d\ge2$, $1\le r\le d$, and $t\ge1$, a rectifiable curve $\gamma:[0,L]\to\mathbb R^d$ with
$\gamma(L)-\gamma(0)\in\mathbb Z^d$, parametrized by arclength,
is called an $r$-sparse $t$-design curve on $\T^d$ if
\[
    \frac1L\int_0^L
        e^{2\pi i\langle m,\gamma(s)\rangle}\,ds
    =
    \int_{\T^d}e^{2\pi i\langle m,x\rangle}\,dx,
    \qquad
\forall\,m\in\Omega_{d,r}(t).
\]
\end{definition}

For fixed $d$ and $r$, we obtain lower bounds of optimal order in $t$ for the length of these curves.
\begin{theorem}\label{thm:sparse-asymptotic}
Let $d\ge2$ and $1\le r\le d$ be integers.
\begin{itemize}
\item[(a)]
If $\gamma$ is an $r$-sparse $t$-design curve on $\T^d$,
for some positive integer $t$, then its length $L$ satisfies
\[
    L\ge
    \left(
        \frac dr\,
        \frac{(t+1)^{2r}-1}{t(t+2)}
    \right)^{1/2}
    \ge t^{r-1}.
\]

\item[(b)]
There is a constant $C_{d,r}>0$ such that, for every positive
integer $t$, there exists a simple $r$-sparse $t$-design
curve on $\T^d$ of length $L\le C_{d,r}t^{r-1}$. 
\end{itemize}
\end{theorem}

\begin{proof}
(a) Since $\|\gamma'(s)\|_2=1$ almost everywhere, the terms 
\[
    E_j:=L\int_0^L|\gamma_j'(s)|^2\,ds,
    \qquad j=1,\ldots,d.
\]
add to $\sum_{j=1}^dE_j=L^2$. 

Fix $I\subseteq\{1,\ldots,d\}$ with $\#I=r$, and define
\[
    A_I:=
    \left\{
        a\in\mathbb R^d:
        a_j\in [t] \text{ for }j\in I,
        \quad a_j=0\ \text{for }j\notin I
    \right\}.
\]
Then $N:=\#A_I=(t+1)^r$ and
$A_I-A_I\subseteq\Omega_{d,r}(t)$.
Consequently, the functions
\[
    f_a(s):=e^{2\pi i\langle a,\gamma(s)\rangle},
    \qquad a\in A_I,
\]
are orthonormal in $L^2([0,L],ds/L)$.

The covariance identity \eqref{eq:11}, applied to the coordinates in $I$, gives
\[
    \frac1N\sum_{a\in A_I}aa^\top
    =
    \frac{t(t+2)}{12}\sum_{j\in I}e_je_j^\top.
\]
Using $f_a'=2\pi i\langle a,\gamma'\rangle f_a$, we obtain
\[
    \frac1N\sum_{a\in A_I}
        \frac1L\int_0^L|f_a'(s)|^2\,ds
    =
    \left(\frac{2\pi}{L}\right)^2
    \frac{t(t+2)}{12}\sum_{j\in I}E_j.
\]

Write $w=\gamma(L)-\gamma(0)\in\mathbb Z^d$.
All the functions $f_a$, $a\in A_I$, satisfy the same
boundary condition,
\[
    f_a(L)=(-1)^{t\sum_{j\in I}w_j}f_a(0).
\]
Thus the periodic or antiperiodic Fourier argument in the proof of Proposition \ref{prop:sobolev} applies to this orthonormal family. Bessel's inequality bounds the total squared Fourier
coefficient at each frequency by one, Parseval's identity gives
total mass $N$, and the sum of the $N$ smallest squared
frequencies is at least $N(N^2-1)/12$, so that we obtain
\[
    \frac1N\sum_{a\in A_I}
        \frac1L\int_0^L|f_a'(s)|^2\,ds
    \ge
    \left(\frac{2\pi}{L}\right)^2
    \frac{N^2-1}{12}.
\]
This implies 
\[
    \sum_{j\in I}E_j
    \ge
    \frac{(t+1)^{2r}-1}{t(t+2)}.
\]
Summing over all $\binom dr$ choices of $I$, and observing
that each coordinate belongs to $\binom{d-1}{r-1}$ of them,
we obtain
\[
    \binom{d-1}{r-1}L^2
    \ge
    \binom dr\,
    \frac{(t+1)^{2r}-1}{t(t+2)}.
\]
The observation $\binom dr/ \binom{d-1}{r-1} = d/r$ yields the first bound in Part (a). The second follows from
\[
    \frac{(t+1)^{2r}-1}{t(t+2)}
    =
    \sum_{\ell=0}^{r-1}(t+1)^{2\ell}
    \ge t^{2r-2}
\]
and $d/r\ge1$.

(b) 
For $r=1$, the primitive vector $v=(1,\ldots,1)^\top$
gives a simple $1$-sparse $t$-design geodesic of length $\sqrt d$.

Suppose $r\ge2$. We use the same counting argument as in the proof of Part (b) in Theorem \ref{thm:hyperbolic asymptotic}. For every integer $1\le s\le t$, counting
frequencies by their support gives
\[
    \#\{m\in\Omega_{d,r}(t):\|m\|_\infty=s\}
    =
    \sum_{h=1}^r\binom dh2^h
        \bigl(s^h-(s-1)^h\bigr)
    \lesssim_{d,r}s^{r-1}.
\]
Consequently,
\[
    \#\Omega_{d,r}(t)\lesssim_{d,r}t^r,
    \qquad
    \sum_{m\in\Omega_{d,r}(t)\setminus\{0\}}
        \frac1{\|m\|_\infty}
    \lesssim_{d,r}
    \sum_{s=1}^t s^{r-2}
    \lesssim_{d,r}t^{r-1}.
\]

Let $\mathcal P$ be the primitive vectors in
$\Omega_{d,r}(t)\setminus\{0\}$, and define
\[
    \mathcal B_k
    :=
    \bigcup_{m\in\mathcal P}
    \{w\in[-k,k]^d\cap\mathbb Z^d:\langle w,m\rangle=0\}.
\]
Dividing a frequency by the greatest common divisor of its
coordinates preserves membership in $\Omega_{d,r}(t)$.
Thus every forbidden hyperplane is represented by an element
of $\mathcal P$.

By the hyperplane-counting bound \eqref{eq:hyperplane-count}
and the union bound,
\[
    \frac{\#\mathcal B_k}{(2k+1)^d}
    \le
    \frac{
        2k\displaystyle\sum_{m\in\mathcal P}\|m\|_\infty^{-1}
        +\#\mathcal P
    }{(2k+1)^2}
    \lesssim_{d,r}
    \frac{t^{r-1}}k+\frac{t^r}{k^2}.
\]
Choose $k=\lceil c_{d,r} t^{r-1}\rceil$. Since $r\ge2$ and $t\ge1$,
\[
    \frac{t^{r-1}}k+\frac{t^r}{k^2}
    \le \frac1c_{d,r}+\frac{t^{2-r}}{c_{d,r}^2}
    \le \frac1c_{d,r}+\frac1{c_{d,r}^2}.
\]
For sufficiently large $c_{d,r}$, the forbidden proportion
is therefore less than one. Hence there exists
$w\in[-k,k]^d\cap\mathbb Z^d$ such that
\[
    \langle w,m\rangle\ne0,
    \qquad
\forall\,m\in\Omega_{d,r}(t)\setminus\{0\}.
\]
The primitive vector $v:=w/\gcd(|w_1|,\ldots,|w_d|)$ 
satisfies the same nonorthogonality conditions and leads to a simple curve of length 
\[
    L:=\|v\|_2
    \le\sqrt d\,k
    \le\sqrt d\,(c_{d,r}+1)t^{r-1}.
\]
Thus, $C_{d,r} = \sqrt d\,(c_{d,r}+1)$ is the required constant.
\end{proof}

For $r=d$, the frequency set $\Omega_{d,r}(t)$ is the full
tensor-product box, and the lower bound in Theorem \ref{thm:sparse-asymptotic} reduces
to Theorem \ref{thm:lower bound}. For $r=1$, the diagonal geodesic has minimum
length $\sqrt d$, independently of $t$. For intermediate
values of $r$, Theorem \ref{thm:sparse-asymptotic} shows that the optimal order of length
is $t^{r-1}$. In particular, the power of \(t\) is determined by the sparsity parameter \(r\), rather than by the ambient dimension \(d\).

\section*{Acknowledgement of AI assistance}
ChatGPT (OpenAI; model: GPT-6 Astra) was used for discussion of selected proof steps and language editing. 
The author takes full responsibility for the content of the paper.

\bibliographystyle{abbrv}

\end{document}